\documentclass{article}
\usepackage{amsmath}
\usepackage{graphicx}
\graphicspath{ {images/} }
\usepackage{bbm}
\usepackage{amsthm}
\usepackage{authblk}
\usepackage{hyperref}
\usepackage{mathrsfs}
\usepackage{amsfonts}
\usepackage{amssymb}
\usepackage{enumerate}
\usepackage[left=3cm,right=2cm,top=2cm,bottom=2cm]{geometry}
\usepackage{epsfig}
 
\usepackage{multicol}
\newcommand*{\email}[1]{%
    \normalsize\href{mailto:#1}{#1}\par
    }
\theoremstyle{definition}

\author{Marcos Morales}
\affil{Pontifical Catholic University of Chile\\ \email{mimorales4@uc.cl}}
\title{Essential Minimum in Families of lines \\}

\begin{document}

\maketitle

\pagestyle{empty}

\setcounter{page}{1}
\pagestyle{plain}
\setlength{\parindent}{0cm}

\begin{abstract}
We apply general methods to generate upper and lower bounds for the essential minimum of a specific family of height functions. In particular, the results shown in this article apply to the case of the Zhang-Zagier height. Furthermore, we can find intervals, where the images of these heights are dense.\\

Our main tool to find upper bounds and intervals of density, is a refinement of the classical Fekete-Szeg\"o theorem due to Burgos Gil, Philippon, Rivera-Letelier and Sombra. \\
\end{abstract}

\begin{center}
\bf{\large{1. Introduction}} \\
\end{center}

In this article we study the  the essential minimum of a family of heights, including the Zhang-Zagier height. We give upper and lower bounds for the essential minimum and intervals of density for the image of these heights.  \\

We now proceed to explain the motivation and describe our results more precisely. To recall the definition of the Zhang-Zagier height, given $\alpha \in \overline{\mathbb{Q}}$, we will consider an algebraic number field  $K$, such that,  $\alpha \in K$. We denote $M_K$ the set of all places of $K$, normalized in such a way that they satisfy the product formula. The Weil height $h: \overline{\mathbb{Q}} \to \mathbb{R}$, is given by

\begin{align*}
h(\alpha) = \displaystyle \frac{1}{[K:\mathbb{Q}]} \sum_{\nu \in M_{K}} \log^+ \left| \alpha \right|_\nu .
\end{align*}

It is well known that $h(\alpha)$ is independent of the field containing $\alpha$ (e.g, see \cite{MR2216774}  Lemma 1.5.2). Then, the Zhang-Zagier height $h_Z : \overline{\mathbb{Q}} \to \mathbb{R}$ is defined by $h_Z(\alpha)= h(\alpha) +h(1-\alpha)$. Our main results concern  the essential minimum of a particular family of curves. Let $f: \overline{\mathbb{Q}} \to \mathbb{R}$ be a function, we define the essential minimum of $f$ as follows 

\begin{center}
$ \mu^{ess}(f)= \inf \lbrace \theta \in \mathbb{R}:\lbrace \alpha \in  \overline{\mathbb{Q}} / f(\alpha) \leq \theta\rbrace $ is an infinite set$ \rbrace $ .  
\end{center} 

It is well known that in the case of the Weil height, $h(\zeta)=0$ for all roots of  unity $\zeta$. Hence, we have  that $\mu^{ess}(h)=0$.\\

Note that, if we take $X \subset \mathbb{C}^2$ a proper subvariety defined over $\mathbb{Q}$, then we can consider the height $h_{X}:X(\overline{\mathbb{Q}}) \to \mathbb{R}$, given by $h_{X}(x,y) = h(x) + h(y)$. The Zhang-Zagier height, corresponds to $V=\lbrace (x,y) \in \mathbb{C}^2 : x+y =1 \rbrace$, then $h_V = h_Z$ . \\

Let $T \subset \mathbb{C}^2$ be a subvariety, we have that $T$ is a torsion subvariety if and only if there are $n,m \in \mathbb{Z} \setminus \left\lbrace 0 \right\rbrace  $ and $\zeta$ a root of unity, such that

\begin{center}
$T= \lbrace (x,y) : x^n y^m = \zeta \rbrace$.
\end{center}

We state a particular case of a Theorem  of S.-W.Zhang.\\

\textbf{\textbf{Theorem.(\cite{MR1189866}):}} \textit{Let $X \subset \mathbb{C}^2$ defined over $\mathbb{Q}$ be a proper subvariety. Then, $\mu^{ess}(h_X)=0$ if and only if $X$ contains a torsion subvariety. }\\

From this theorem we can conclude that $\mu^{ess}(h_Z)>0$. The question now, is how to compute the essential minimum of heights like $h_Z$. However, the answer remains unknown. Despite of this, there have been several attempts to approach this number. In $1993$ Zagier proved the following theorem \\

\textbf{\textbf{Theorem.(\cite{MR1197513}):}} \textit{ For all $\alpha \in \overline{\mathbb{Q}}$ such that $\alpha \notin \left\lbrace 0,1, e^{i\pi/3},e^{-i\pi/3}\right\rbrace $, we have}

\begin{center}
$h_Z(\alpha) \geq \displaystyle \frac{1}{2}\log \left( \frac{1+\sqrt{5}}{2}\right) \approx 0.2406059...$.
\end{center}

This theorem tells us that $\mu^{ess}(h_Z) \geq 0.240606$. In 2001  Doche improved this result and proved the following theorem  \\

\textbf{\textbf{Theorem.(\cite{MR1838073}):}} \textit{ Let $h_Z$ be the Zhang-Zagier height, then}

\begin{center}
$ 0.2482474 \leq \mu^{ess}(h_Z) \leq 0.25443678$.
\end{center}

This result is the best approximation for $\mu^{ess}(h_Z)$ that has been computed until now. \\ 

Our contribution does not improve Doche's theorem, however, it gives an interval where the image of the height $h_Z$ is dense. \\

\textbf{\textbf{Theorem A:}} \textit{ Let $h_Z$ be the Zhang-Zagier height. Then, the image of $h_Z$ is dense in the interval $[0.31944909, \infty)$. In particular}

\begin{center}
$ \mu^{ess}(h_Z) \leq 0.31944909$.
\end{center}

Theorem A is a direct consequence of Theorem 1 in  \cite{riveraspectrum}. Also, Theorem A is consequence of a more general result, Theorem B below. For $a,b \in \overline{\mathbb{Q}}$ we consider the subvariety $L_{a,b} = \lbrace (x,y) \in \mathbb{C}^2 : y= ax +b$ and $a,b \in \overline{\mathbb{Q}}\rbrace$, and we denote $h_{a,b} = h_{L_{a,b}}$. We compute upper and lower bounds for $\mu^{ess}(h_{a,b})$. Before proceeding, we introduce some notations. Let $a,b \in \overline{\mathbb{Q}}$, we write $K_a$ the Galois closure of $\mathbb{Q}(a)/\mathbb{Q}$ and $K_{a,b}$ the Galois closure of the field generated by $a$ and $b$ over $\mathbb{Q}$. We also write $G(a)=$Gal$(K_a/\mathbb{Q})$, $G(a,b)=$Gal$(K_{a,b}/\mathbb{Q})$, $\deg(a)= [\mathbb{Q}(a): \mathbb{Q}]$ and Gal$(a)= \lbrace \sigma(\alpha) : \sigma \in G(a)\rbrace$. \\

Let $\sigma \in  G(a,b)$, we define,  $\Psi_{a,b}^\sigma : \mathbb{R} \to \mathbb{R}$, $\varphi: \mathbb{R} \to \mathbb{R}$, $\Delta: \overline{\mathbb{Q}} \times \overline{\mathbb{Q}} \to \mathbb{R}$ and  $\Omega_{a,b}: \mathbb{R} \to \mathbb{R}$ , given by

\begin{equation}
\Psi_{a,b}^\sigma(t) =  \displaystyle \frac{1}{2\pi} \int_0^{2\pi}\log^{+} \left| e^{i\theta} + \frac{\sigma(b)}{\sigma(a)} + t\right| d\theta,
\end{equation}
\begin{equation}
\varphi(t) =\displaystyle \frac{1}{2\pi} \int_{0}^{2 \pi} \log^{+} |e^{i\theta} + t|d\theta,
\end{equation}
\begin{equation}
\Delta(a,b) =\displaystyle \sum_{\substack{ p \\ \text{prime}}} \sum_{\sigma \in G(a,b)} \log^+  \max(|\sigma(a)|_p,|\sigma( b)|_p) +\displaystyle \sum_{\sigma \in G(a,b)} \log^+ |\sigma(a)| ,
\end{equation}

\begin{equation}
\Omega_{a,b}(t) = \Delta(a,b)+\varphi(t) \hspace{1mm}+ \displaystyle  \sum_{\sigma \in G(a,b)} \Psi_{a,b}^\sigma(t).\\
\end{equation}

Now, we can state the following theorem \\

\textbf{\textbf{Theorem B:}} \textit{ Let $a,b \in \overline{\mathbb{Q}}$. Then, there exists an effectively computable number $\mathcal{K}(a,b)$  such that for each $t \in \mathbb{R}$, we have that}

\begin{center}
$ \mathcal{K}(a,b) \leq \mu^{ess}(h_{a,b}) \leq \Omega_{a,b}\left(  \displaystyle t \right) $.
\end{center}

\textit{Furthermore, assume that $a$ and $b$ satisfy one of the following properties}\\

\begin{enumerate}[i)]
\item  \textit{There exists $\sigma_0\in G(a,b)$, such that, $||\sigma_0(a)|-|\sigma_0(b)|| > 1$.} 
\item  \textit{ There exists $\sigma_0\in G(a,b)$, such that, $||\sigma_0(a)|-|\sigma_0(b)|| < 1$ and $|\sigma_0(a)|+|\sigma_0(b)|< 1$.} 
\end{enumerate}

\textit{Then, $\mathcal{K}(a,b) > 0$. See section $3$ for the definition of $\mathcal{K}(a,b)$} \\

The upper bound obtained in Theorem $A$ is a particular case of theorem $B$, with $a=-1$ and $b=1$. In this case, the best upper bound found uses $t=-1/2$. \\ 

 Let $a,b \in \mathbb{Q}$, $a\neq 0$. If $b \neq 0$, we can write $a= a_1/a_2$, $b=b_1/b_2$  and $(a_1,a_2)=1$, $(b_1,b_2)=1$. We define $S_{a,b}= \lbrace p $ prime $: p|a_2 \vee p|b_2\rbrace$. If $b=0$, then $S_{a,0}= \lbrace p : p|a_2 \rbrace$. We also write $s= |S_{a,b}|$, then $S=  \lbrace p_i : 1 \leq i \leq s \rbrace$. Using this notation, we can define the function $\Gamma_{a,b}: \mathbb{R}^{s+1} \to \mathbb{R}$

\begin{align*}
\Gamma_{a,b} (x,r_1,r_2,...,r_s)= \displaystyle \sum_{i=1}^s \log^+ |r_i|+\log^+ \max (|a|_{p_i} r_i, |b|_{p_i}) + \frac{1}{2\pi} \int_0^{2\pi} \log^+ \left| \frac{e^{i\theta}}{r_1r_2... r_s} +x \right| + \log^+\left|\frac{ae^{i\theta}}{r_1r_2... r_s}+b+x \right| d\theta.
\end{align*}

\textbf{\textbf{Theorem C:}} \textit{ Let $a,b \in \mathbb{Q}$. Then, for each $x,r_1,r_2,..., r_s \in \mathbb{R}$, the image of $h_{a,b}$ is dense in the interval $[\Gamma_{a,b}(x,r_1,r_2,..., r_s), \infty)$, in particular}
\begin{align*}
\mu^{ess} (h_{a,b}) \leq \Gamma_{a,b}(x,r_1,r_2,..., r_s).
\end{align*}

This last theorem does not improve the upper bound given in Theorem $B$, however, it does give intervals of density for $a,b \in \mathbb{Q}$. For instance, if we take $a=1$ and $b=2$, we find that $\mathcal{K}(1,2)= \log(\sqrt{3})$ and $\Omega_{1,2}(0) \leq 0.6461599$ (see Theorem 2.9 and Corollary 3.9). Therefore, Theorem B and Theorem C, give us that $ \log(\sqrt{3})= 0.5493061... \leq \mu^{ess}(h_{1,2}) \leq 0.6461599$ and the image of $h_{1,2}$, is dense in the interval $[0.6461599, \infty]$.\\ \\ 

We will start in section $2$ by determining upper bounds for the essential minimum. Then, in section $3$ we will compute lower bounds and prove Theorem B using the methods outlined in \cite{MR3802441}. Finally, in section $4$ we will determine intervals of density for the image of $h_{a,b}$ for $a,b \in \mathbb{Q}$  and prove Theorem A and Theorem C using results from other sections. \\ \\

\textbf{Acknowledgements} \\

I want to thank to my Ph.D. advisor Ricardo Menares and the maths department of Pontificia Universidad Católica de Chile. \\ \\

\begin{center}
\textit{\bf{2. Upper Bounds}}
\end{center}

Our purpose is to give a good upper bound for the essential minimum of each element in the family of heights $\left\lbrace h_{a,b} \right\rbrace_{a,b \in \overline{\mathbb{Q}} } $. Since the case $a=0$ is trivial, we will assume henceforth $a \neq 0$.\\

Let $p$ be a prime and let $|.|_p$ be the standard $p-$adic value on $\mathbb{Q}_p$. It is well known that $|.|_p$ can be uniquely extended to $\overline{\mathbb{Q}}_p$. We fix and embedding $\iota_0 : \overline{\mathbb{Q}} \to \overline{\mathbb{Q}}_p$, and for $\alpha \in \overline{\mathbb{Q}}$, we define $|\alpha|_p = | \iota_0(\alpha)|_p$. Before proceeding, we need the following  lemmas. \\

\textit{{\bf{Lemma 2.1:}} Let $\alpha \in \overline{\mathbb{Q}}$, and $K$  an algebraic number field containing $\alpha$. Then} \\

\begin{align*}
h(\alpha) = \displaystyle \frac{1}{[K:\mathbb{Q}]}  \left( \sum_{\substack{ p \\ \text{prime}}}\sum_{ \sigma \in \text{Gal}(K/\mathbb{Q})} \log^+ \left| \sigma(\alpha) \right|_p  + \sum_{\sigma \in \text{Gal}(K/\mathbb{Q})} \log^+ \left|\sigma(\alpha) \right|_{\infty} \right).
\end{align*}

\textit{Proof: } Using \cite{MR2216774}, Corollary 1.3.5 and its proof, we have that, for each place  $\nu \in M_K$, there exists a unique $\sigma \in$ Gal$(K/\mathbb{Q})$ and a unique place $|.|_w$ of $\mathbb{Q}$, such that $|.|_\nu = |. |_w \circ \sigma$, this proves the lemma.$\hfill\square$ \\ \\

Now, we define the function $U_{a}^b: \overline{\mathbb{Q}} \to \mathbb{R} $, given by 

\begin{align*}
U_{a}^b(\alpha) = \displaystyle \frac{1}{\deg(\alpha)} \displaystyle \sum_{\substack{\beta \in \text{Gal}(\alpha) \\ \sigma \in G(a,b)}}  \log^+ \left|\beta + \frac{\sigma(b)}{\sigma(a)}\right|.
\end{align*}

\textit{{\bf{Lemma 2.2:}} Let $a,b \in \overline{\mathbb{Q}}$ and $\alpha \in \overline{\mathbb{Z}}$, then}

\begin{center}
$h(a\alpha +b) \leq U_a^b(\alpha)+  \Delta(a,b).$
\end{center}

\textit{Proof: }We consider $K$ the Galois closure of the field generated by $a$, $b$ and $\alpha$. It is clear that $K_a \subseteq K_{a,b} \subset K$ and $\deg(\alpha) \leq  [K_{\alpha}:\mathbb{Q}] \leq [K:\mathbb{Q}]$. Then, using Lemma 2.1, we have that
\begin{align*}
h(a\alpha +b) &= \displaystyle \frac{1}{[K:\mathbb{Q}]} \left( \sum_{\substack{p \\ \text{prime}}}\sum_{ \sigma \in \text{Gal}(K/\mathbb{Q}) } \log^+ \left| \sigma (a\alpha +b) \right|_p  + \sum_{\sigma \in \text{Gal}(K/\mathbb{Q})} \log^+ \left| \sigma(a\alpha +b) \right| \right)\\
&\leq \displaystyle \frac{1}{[K_{\alpha}:\mathbb{Q}]} \left( \sum_{\substack{p \\ \text{prime}}}\sum_{\substack{\delta \in G(\alpha) \\ \sigma \in G(a,b)}} \log^+ \left| \sigma(a)\delta(\alpha) +\sigma(b) \right|_p  + \sum_{\substack{\delta \in G(\alpha) \\ \sigma \in G(a,b)}}  \log^+ \left|\sigma(a)\delta(\alpha) + \sigma(b)\right| \right). \\
\end{align*}

In the last inequality we have used that $[K_{\alpha}:\mathbb{Q}] \leq [K:\mathbb{Q}]$ and that the number of elements in $\text{Gal}(K/\mathbb{Q})$ is less than or equal to the number of pairs $(\sigma, \delta)$, with $\sigma \in G(\alpha)$ and $\delta \in G(a,b)$. Using the fact that $\alpha \in \overline{\mathbb{Z}}$, we conclude that for every $p$ prime, $\sigma \in G(a,b)$ and $\delta \in G(\alpha)$, we have $|\sigma(a)\delta(\alpha) + \sigma(b)|_p \leq \max(|\sigma(a)|_p,|\sigma(b)|_p)$. Therefore

\begin{align*}
h(a\alpha +b) &\leq \displaystyle \sum_{\substack{p \\ \text{prime}}} \sum_{\sigma \in G(a,b)} \log^+  \max(|\sigma(a)|_p,|\sigma(b)|_p) +  \frac{1}{[K_{\alpha}:\mathbb{Q}]} \left( \displaystyle \sum_{\substack{\delta \in G(\alpha) \\ \sigma \in G(a,b)}}  \log^+ \left|\sigma(a)\right| + \log^+ \left|\delta(\alpha) + \frac{\sigma(b)}{\sigma(a)}\right| \right)\\
&\leq  \displaystyle \sum_{\substack{p \\ \text{prime}}} \sum_{\sigma \in G(a,b)} \log^+  \max(|\sigma(a)|_p,|\sigma(b)|_p) + \displaystyle \sum_{\sigma \in G(a,b)} \log^+ |\sigma(a)| + \frac{1}{[K_{\alpha}:\mathbb{Q}]} \displaystyle \sum_{\substack{\delta \in G(\alpha) \\ \sigma \in G(a,b)}}  \log^+ \left|\delta(\alpha) + \frac{\sigma(b)}{\sigma(a)}\right| \\
& =  \displaystyle \sum_{\substack{p \\ \text{prime}}} \sum_{\sigma \in G(a,b)} \log^+  \max(|\sigma(a)|_p,|\sigma(b)|_p) + \displaystyle \sum_{\sigma \in G(a,b)} \log^+ |\sigma(a)| + \frac{1}{\deg(\alpha)} \displaystyle \sum_{\substack{\beta \in \text{Gal}(\alpha) \\ \sigma \in G(a,b)}}  \log^+ \left|\beta + \frac{\sigma(b)}{\sigma(a)}\right| \\
&= U_a^b(\alpha)+  \Delta(a,b).
\end{align*}
 
This concludes the proof of the lemma.$\hfill\square$ \\ \\

Let $E \subset \mathbb{C}$, we denote $\text{Cap}(E)$, the capacity of $E$ and $\mu_E$ the equilibrium measure of $E$, see \cite{MR1009368} section $3$, for the definitions. Let $z \in \mathbb{C}$, we define $d(z,E)= \inf_{a \in E} |z-a|$, and $B(E,r)= \lbrace z \in \mathbb{C}: d(z,E) < r\rbrace$ for $r >0$.\\

Given a compact set $E \subset \mathbb{C}$, we denote $C_0(E, \mathbb{R})$, the set of all continuous functions from $E$ to $\mathbb{R}$. Let $E= \lbrace a_1,a_2,...,a_k \rbrace \subset \mathbb{C}$ a finite set, we define the measure $\delta(E): 	C_0(E,\mathbb{R}) \to \mathbb{R}$, by

\begin{center}
$\delta(E)(f) = \displaystyle \frac{1}{k}  \sum_{n=1}^k f(a_n)$.
\end{center}

Now, we can set the following proposition\\

\textit{{\bf{Proposition 2.3: }}Let $E \subset \mathbb{C}$ a compact set invariant under complex conjugation, then, there exists a sequence of algebraic integers $\lbrace \alpha_n \rbrace_{n \in \mathbb{N}}$, such that $\text{Gal}(\alpha_n) \subset B\left( E, \frac{1}{n}\right)$, and}\\ 

\begin{center}
$\delta(Gal(\alpha_n)) \overset{*}{\longrightarrow} \mu_E$.
\end{center}

\textit{Proof: }Using \cite{MR3928039}, Proposition 7.4, with $E_{|.|_\infty}=E$ and $E_{|.|_p}= \mathcal{O}_p= \lbrace z \in \mathbb{C}_p : |z|_p \leq 1 \rbrace$, we conclude the proof. $\hfill\square$\\

\textit{{\bf{Proposition 2.4:}} Let $a,b \in \overline{\mathbb{Q}}$, and $a \neq 0$. Then, for each $t \in \mathbb{R}$ we have that}
\begin{align*}
\mu^{ess} (h_{a,b}) \leq \Omega_{a,b}(t).
\end{align*}

\textit{Here, $\Omega_{a,b}$ is the function defined in $(4)$.}\\
\hspace{1cm}

\textit{Proof:} By Lemma 2.2, given $\alpha \in \overline{\mathbb{Z}} $, we have that

\begin{align*}
h_{a,b}(\alpha) &=  h(\alpha) + h(a\alpha + b), \\
&\leq h(\alpha) + U_a^b(\alpha) +  \Delta(a,b) =: \eta_{a,b}(\alpha).
\end{align*}

For Proposition $2.3$, given $E \subseteq \mathbb{C}$, a compact set with $\text{Cap}(E) =1$ and invariant under complex conjugation, there exists a sequence of algebraic integers $\alpha_n \in \overline{\mathbb{Z}}$, such that $\text{Gal}(\alpha_n) \subset B\left( E, \frac{1}{n}\right)$ and

\begin{align*}
h(\alpha_n) + U_a^b(\alpha_n)  \underset{n \to \infty}{\longrightarrow} \displaystyle\int_{E} \log^{+}\left | x \right | d \mu_E(x)  +\displaystyle \sum_{\sigma \in G(a,b)} \int_{E} \log^{+}\left | x +\frac{\sigma(b)}{\sigma(a)} \right | d \mu_E(x) =: M_E.
\end{align*}

Therefore

\begin{align*}
h_{a,b}(\alpha_n) \leq \eta_{a,b}(\alpha_n) \underset{n \to \infty}{\longrightarrow}  \Delta(a,b) +M_E =: J_E.
\end{align*}

We claim that $J_E < \infty$. In fact, we have that the only possible unbounded term in the definition could be $M_E$. Let the functions, $f: \mathbb{R} \to \mathbb{R}$, and for each $\sigma \in G(a,b)$ the function $g_{\sigma}: \mathbb{R} \to \mathbb{R}$, be defined by  $f(t)=\log^+|t|$ and $g_\sigma(t)=\log^+| t+ \sigma(b)/\sigma(a)|$. These functions are continuous and $E$ is a compact set, therefore, $M_E < \infty$. Since $J_E < \infty$, the sequence $\lbrace h(\alpha_n)\rbrace_{n \in \mathbb{N}}$ is bounded, we conclude that there is a subsequence which is convergent, we call it $\lbrace\beta_n\rbrace_{n \in \mathbb{N}}$. \\

If $h_{a,b}(\beta_n) \underset{n \to \infty}{\longrightarrow} Z$, then, by definition of limit, given $\varepsilon > 0$ the set $\lbrace \beta_n \in \overline{\mathbb{Z}} : h_{a,b}(\beta_n) \leq Z +\varepsilon \rbrace$ is infinite, therefore for every $\varepsilon > 0$ we have $\mu^{ess}(h_{a,b}) \leq Z + \varepsilon$. Taking $\varepsilon \to 0$ we get $\mu^{ess}(h_{a,b}) \leq Z$. Since $Z \leq J_E$, we conclude that $\mu^{ess}(h_{a,b}) \leq J_E$. Given  $t \in \mathbb{R}$ we use $E=S_t= S_1 +t$, where $S_1 = \lbrace z \in \mathbb{C} : |z|=1 \rbrace$. Then, $\mu_{S_t}$ is the natural translation of the measure $\mu_{S_1}=\frac{d\theta}{2\pi}$. We deduce that

\begin{center}

$\mu^{ess}(h_{a,b}) \leq J_{S_t} = \Omega_{a,b}(t)$.

\end{center}

This concludes the proof of the theorem $\hfill\square$ \\ \\

Let $a,b \in \mathbb{Q}$, $a\neq 0$, then

\begin{equation}
\Omega_{a,b}(t) = \Delta(a,b) + \varphi(t) +\varphi(t +b/a).\\
\end{equation}

We want to find a point for which $\Omega_{a,b}$ achieves its minimum value and compute a power series for $\Omega_{a,b}$ at that point. Before proceeding, we will prove the following two lemmas  \\

\textit{{\bf{Lemma 2.5:}} The function $\varphi$ satisfies the following properties} 

\begin{enumerate}[(a)]
\item[i)] \textit{For each $t \in \mathbb{R}$, $\varphi(t) = \varphi (-t)$,}
\item[ii)] \textit{for $ |t| \geq 2$, $\varphi(t) = \log |t|$.}
\end{enumerate}

\textit{Proof:} Let $t \in  \mathbb{R}$, we have

\begin{align*}
\varphi(-t) &=   \displaystyle \frac{1}{2\pi} \int_0^{2\pi} \log^{+} \left| e^{i\theta} - t\right| d\theta
 = \displaystyle \frac{1}{2\pi} \int_0^{2\pi} \log^{+} \left| e^{i(\theta +\pi)} +t\right| d\theta\\ 
 &= \displaystyle \frac{1}{2\pi} \int_{\pi}^{3\pi} \log^{+} \left| e^{i\theta} +t\right| d\theta 
  = \displaystyle \frac{1}{2\pi} \left(  \int_{\pi}^{2\pi} \log^{+} \left| e^{i\theta} +t\right| d\theta +\int_{2\pi}^{3\pi} \log^{+} \left| e^{i\theta} +t\right| d\theta \right)\\
  &= \displaystyle \frac{1}{2\pi} \left(  \int_{\pi}^{2\pi} \log^{+} \left| e^{i\theta} +t\right| d\theta +\int_{0}^{\pi} \log^{+} \left| e^{i\theta} +t\right| d\theta \right) \\
    &= \displaystyle \frac{1}{2\pi} \int_{0}^{2\pi} \log^{+} \left| e^{i\theta} +t\right| d\theta = \varphi(t). \hspace{7cm}  
\end{align*}

This proves (i). Now, if $|t| \geq 2$, then for each $\theta \in \mathbb{R}$, $|e^{i\theta} +t| > |1-|t|| = |t|-1 \geq 1$, hence, $\log^+ |e^{i\theta} +t| = \log |e^{i\theta}+t|$, we conclude that

\begin{align*}
\varphi(-t)&= \displaystyle\frac{1}{2\pi}  \displaystyle\int_{0}^{2\pi} \log\left |t-e^{i\theta}\right | d\theta\\
&= \displaystyle\frac{1}{2\pi} Re \left( \displaystyle\int_{0}^{2\pi} \log\left (t-e^{i\theta}\right ) d\theta \right)\\
&=  \displaystyle\frac{1}{2\pi} Re \left( \displaystyle\int_{S_1} \frac{1}{iz}\log\left (t-z\right ) dz \right)\\
&=  \displaystyle\frac{1}{2\pi} Re \left( 2\pi i\frac{1}{i}\log\left (t\right )\right)\\
&=  \log |t|\\
\end{align*}
This proves (ii).  $\hfill\square$ \\ \\

\textit{{\bf{Lemma 2.6:}} Let $f: \overline{\mathbb{Q}} \to \mathbb{R}$, $\gamma, \delta \in \overline{\mathbb{Q}}$ and $\gamma \neq 0$, define $f^{\gamma,\delta} :\overline{\mathbb{Q}} \to \mathbb{R}$ by $f^{\gamma,\delta} (\alpha) =f(\gamma \alpha +\delta)$, then }

\begin{center}
$\mu^{ess}(f) =\mu^{ess}(f^{\gamma,\delta})$.
\end{center}

\textit{Proof:}  In fact, if $\mu^{ess}(f) = A$, then there exist a sequence ${\alpha_n}_{n\in \mathbb{N}} \subseteq \overline{\mathbb{Q}}$ such that $f(\alpha_n) \to A$. Let $\beta_n = (\alpha_n-\delta)/\gamma \in \overline{\mathbb{Q}}$, then $f^{\gamma,\delta} (\beta_n) = f(\alpha_n) \to A $, therefore $\mu^{ess}(f^{\gamma,\delta}) \leq \mu^{ess}(f)$. Using the fact that $1/\gamma, -\delta \in \overline{\mathbb{Q}}^*$ and the same argument, we can show the other inequality and the lemma is proved. $\hfill\square$\\

\textit{{\bf{Corollary 2.7:}} Let $a,b \in \mathbb{Q}$, then }

\begin{center}
$\mu^{ess}(h_{a,b})= \mu^{ess}(h_{|a|,|b|})$.
\end{center}

\textit{Proof:}  We have that $h_{-a,b} (\alpha) =h(\alpha) + h(-a\alpha + b)= h(-\alpha) + h(a(-\alpha) + b) = h_{a,b}(-\alpha)$. Using Lemma 2.6 with $\gamma= -1$ and $\delta=0$, we conclude that $\mu^{ess}(h_{-a,b}) = \mu^{ess}(h_{a,b})$. We also note that $h_{a,-b}(\alpha)=h(\alpha) + h(a\alpha - b)= h(-\alpha) + h(-a\alpha + b) = h_{-a,b}(\alpha) = h_{a,b}(-\alpha) $. Hence, we have $\mu^{ess}(h_{-a,b}) = \mu^{ess}(h_{a,-b}) = \mu^{ess}(h_{a,b}) $. Combining this two equalities we can conclude all the others, this completes the proof of the corollary.  $\hfill\square$ \\ \\

Let $c \in \mathbb{R}$, we define the function $\varepsilon_{c} (z,n): \mathbb{N} \times \mathbb{R} \to \mathbb{C}$, given by

\begin{center}
$\varepsilon_{c}(n,t) = \displaystyle\int_{c}^{\pi} \left( e^{i\theta}-e^{it}\right)^n d\theta$ ;
\end{center}

Expanding the binomial and integrating, we get

\begin{center}
$ \varepsilon_c (n,t) = (-1)^n(\pi-c)e^{itn} + \displaystyle\sum_{k=0}^{n-1} \frac{1}{i(n-k)} {n \choose k} (-1)^k(e^{i\pi(n-k)+itk}-e^{ic(n-k)+itk})$ .
\end{center}

Using this notation, we have the following proposition. \\

\textit{{\bf{Proposition 2.8:}}} Let $a,b \in \mathbb{Q}$ with $a > 0$ and $b \geq 0$. 

\begin{enumerate}[(a)]
\item[i)] \textit {If $b=0$, then}
\begin{center}
$\Omega_{a,0}(0)= h(a)$.
\end{center}
\item[ii)] \textit{If $0 < b/a < 4$, then}
\begin{align*}
\Omega_{a,b}\left( \displaystyle -\frac{b}{2a}\right) = \Delta(a,b) + \displaystyle\frac{2}{\pi}  Re\left(   \log \left( \frac{b}{2a} - e^{i\frac{\pi + \alpha_{a,b}}{2}}\right)(\pi-\alpha_{a,b}) - \displaystyle\sum_{n=1}^{\infty} \frac{(2a)^n}{n \left(b-2ae^{i\frac{\pi + \alpha_{a,b} }{2}}\right)^n} \varepsilon_{\alpha_{a,b}} \left( n, \frac{\pi + \alpha_{a,b}}{2} \right) \right).
\end{align*} 

\textit {Where $\alpha_{a,b}= \arctan\left( \frac{\sqrt{16a^2-b^2}}{b}\right)$ and $\log$ is the main branch of the logarithm.}

\item[iii)] \textit{ If $b/a \geq 4$, then}
\begin{center}
$\Omega_{a,b} \displaystyle  \left(  0 \right) = \Delta(a,b) + \displaystyle \log\left( \frac{b}{a}\right)$.
\end{center}
\end{enumerate}

\textit{Proof:} If $b=0$, using (5) and $\varphi(0)$=0, we get  that $\Omega_{a,b}(0)= h(a)$. Now, assume $ 0< b/a <4$. Using (4) and Lemma 2.5 (i), we have

\begin{equation}
\displaystyle\Omega_{a,b}\left( -\frac{b}{2a}\right) = \Delta(a,b) +2 \varphi\left(-\frac{b}{2a}\right)  = \Delta(a,b) +\displaystyle\frac{2}{\pi}  \displaystyle\int_{\alpha_{a,b}}^{\pi} \log\left |\frac{b}{2a}-e^{i\theta}\right | d\theta.
\end{equation}

Where $\alpha_{a,b}$ is the argument of the complex number given by the intersection of the two circumferences $S_{\frac{b}{2a}}$ and $S_1$ in the first quadrant. \\

Let $f: \mathbb{C} \setminus \left[ \frac{b}{2a},+\infty \right)  \rightarrow \mathbb{C}$, defined by $f(z)= \log \left( \frac{b}{2a}-z \right)$, where $\log$ is the main branch of the logarithm, then  the power series of $f$ around $z=e^{i\frac{\pi + \alpha_{a,b}}{2}}$, is given by

\begin{center}
$f(z)= \displaystyle\log \left(\frac{b}{2a} - e^{i\frac{\pi + \alpha_{a,b}}{2}}\right) + \displaystyle\sum_{n=1}^{\infty} \frac{-(2a)^n}{n \left(b-2ae^{i\frac{\pi + \alpha_{a,b}}{2}}\right)^n}\left( z-e^{i\frac{\pi + \alpha_{a,b}}{2}} \right) ^n$.
\end{center}

The convergence radius of this series is 

\begin{align*}
\displaystyle\frac{1}{r} &= \displaystyle\underset{n \to \infty}{\displaystyle\lim} \displaystyle \left |  \frac{2na}{\left( n+1\right) \left( b-2ae^{i\frac{\pi + \alpha_{a,b}}{2}} \right) }  \right |= \frac{2a}{\left | b-2ae^{i\frac{\pi + \alpha_{a,b}}{2}} \right |} ,\\
\displaystyle r &= \displaystyle\left | \frac{b}{2a}-e^{i\frac{\pi + \alpha_{a,b}}{2}} \right |.
\end{align*}

On the other hand, we have

\begin{align*}
\displaystyle\int_{\alpha_{a,b}}^{\pi} \log \left |\frac{b}{2a}-e^{i\theta} \right | d\theta &= Re\left( \displaystyle\int_{\alpha_{a,b}}^{\pi} \log\left( \frac{b}{2a}-e^{i\theta}\right)  d\theta\right)\\
 &= Re\left( \displaystyle\int_{\alpha_{a,b}}^{\pi} \log\left( \frac{b}{2a} - e^{i\frac{\pi + \alpha_{a,b}}{2}}\right)  + \displaystyle\sum_{n=1}^{\infty} \frac{-(2a)^n}{n \left(  b-2ae^{i\frac{\pi + \alpha_b}{2}}\right) ^n}\left( e^{i\theta}-e^{i\frac{\pi + \alpha_{a,b}}{2}}\right)^n d\theta\right).
\end{align*}

The maximum value of $\left | e^{i\theta}-e^{i\frac{\pi + \alpha_{a,b}}{2}} \right|$ for $\theta \in [\alpha_{a,b},\pi]$  is achieved when $\theta= \pi$, therefore for  $\theta \in [\alpha_{a,b},\pi]$, we have

\begin{center}
$ \left | e^{i\theta}-e^{i\frac{\pi + \alpha_{a,b}}{2}} \right| \leq \left | 1+e^{i\frac{\pi + \alpha_{a,b}}{2}} \right|$ .
\end{center}

Note that $0< \alpha_{a,b} \leq \pi/2$, therefore $\pi/2 < (\pi + \alpha_{a,b})/2 \leq 3\pi/4$, it follows that   \\

\begin{center}
$  -\displaystyle\frac{\sqrt{2}}{2} \leq \displaystyle \cos \left( \frac{\pi + \alpha_{a,b}}{2} \right) < 0.$
\end{center}

From this last equation, we conclude that

\begin{equation}
 \left | e^{i\theta}-e^{i\frac{\pi + \alpha_{a,b}}{2}} \right| \leq \left | 1+e^{i\frac{\pi + \alpha_{a,b}}{2}} \right| < \displaystyle\left | \frac{b}{2a}-e^{i\frac{\pi + \alpha_{a,b}}{2}} \right | = r.
\end{equation}

Since the convergence is uniform we can exchange the integral with the series in $(5)$ and we get

\begin{align*}
\Omega_{a,b}\left(  \displaystyle-\frac{b}{2a}\right)  = \Delta(a,b)+\displaystyle\frac{2}{\pi} Re\left(   \log \left( \frac{b}{2a} - e^{i\frac{\pi + \alpha_{a,b}}{2}}\right)(\pi-\alpha_{a,b}) -  \displaystyle\sum_{n=1}^{\infty} \frac{(2a)^n}{n \left(b-2ae^{i\frac{\pi + \alpha_{a,b} }{2}}\right)^n} \varepsilon_{\alpha_{a,b}} \left( n , \frac{\pi +\alpha_{a,b}}{2}\right)  \right).
\end{align*}

On the other hand, if $|b/a| \geq 4$, then $|b/a| > 2$, using (ii) from lemma 2.4, we conclude that, $\Omega_{a,b}(0) = \Delta(a,b) + \varphi(0)+ \varphi(b/a)= \Delta(a,b) + \log (b/a)$,  this completes the proof of the proposition. $\hfill\square$ \\ \\

Using Proposition 2.4 and Proposition 2.9, we can prove the following theorem\\

\textit{{\bf{Theorem 2.9:}} Let $a,b \in \mathbb{Q}$. If $b=0$, we have that }

\begin{center}
$\mu^{ess}(h_{a, 0}) \leq \Omega_{|a|, 0} \left( 0 \right)  = h(a)$.
\end{center} 

For $|b/a|=1$, we have that

\begin{center}
$\mu^{ess}(h_{a, b}) \leq \Omega_{|a|, |b|} \left( -\frac{1}{2} \right)  \leq \Delta(a,b)+ 0.3194490869562$.

\end{center} 

For $|b/a|=2$, 
\begin{center}
$\mu^{ess}(h_{a, b}) \leq \Omega_{|a|, |b|} \left( -1 \right)  \leq  \Delta(a,b) +0.6461598436469$.

\end{center} 

For $|b/a|=3$, 
\begin{center}
$\mu^{ess}(h_{a,b}) \leq \Omega_{|a|,|b|} \left( - \frac{3}{2} \right) \leq \Delta(a,b) +0.9909205628144$.

\end{center}

For $|b/a| \geq 4$,

\begin{center}
$\mu^{ess}(h_{a,b}) \leq \Omega_{|a|, |b|} \left( 0 \right) = \Delta(a,b) + \displaystyle  \log\left(\frac{  b  }{a} \right) $.
\end{center}

\textit{Proof:} Using Corollary $2.7$, we can assume that $a>0$ and $b \geq 0$. Assume first that, $b/a=1$. Let

\begin{align*}
T_1 &:= \displaystyle\frac{2}{\pi} Re\left( \log \left( \frac{1}{2} - e^{i\frac{\pi + \alpha_{a,b}}{2}}\right)(\pi-\alpha_{a,b}) -  \displaystyle\sum_{n=1}^{20} \frac{2^n}{n \left(1-2e^{i\frac{\pi + \alpha_{a,b} }{2}}\right)^n} \varepsilon_{\alpha_{a,b}} \left( n, \frac{\pi + \alpha_{a,b}}{2} \right) \right)\\
 &= 0.3194345111561... \leq  0.3194345111562 .
\end{align*}
 
Let $R_1=   \Omega_{a,b}\left( -\frac{1}{2}\right) - T_1 - \Delta(a,b)$. Using Proposition 2.8, we have that

\begin{align*}
|R_1| &\leq \displaystyle\frac{2}{\pi}\displaystyle\sum_{n=21}^{\infty} \frac{2^n}{n \left|1-2e^{i\frac{\pi + \alpha_{a,b}}{2}}\right|^n}\displaystyle\int_{\alpha_{a,b}}^{\pi}\left| e^{i\theta}-e^{i\frac{\pi + \alpha_{a,b}}{2}}\right| ^n d\theta \\
\end{align*}

Note that, since $b/a= 1$, we have $\alpha_{a,b}= \arctan(\sqrt{15})$, therefore $ \left|1-2e^{i\frac{\pi + \alpha_{a,b}}{2}}\right| = \sqrt{5 +\sqrt{6}}$. Moreover, using $(7)$, we conclude that, for each $\theta \in [\alpha_{a,b}, \pi]$, $\left| e^{i\theta}-e^{i\frac{\pi + \alpha_{a,b}}{2}}\right| \leq \sqrt[]{2-\sqrt[]{\frac{3}{2}}} $, therefore

\begin{align*}
|R_1| &\leq \displaystyle\frac{2}{\pi}\displaystyle\sum_{n=21}^{\infty} \frac{2^n}{n \left( \sqrt{5+\sqrt{6}}\right)  ^n}\left( \sqrt[]{2-\sqrt[]{\frac{3}{2}}}\right) ^n (\pi-\alpha_{a,b}) \\
 &\leq \displaystyle\frac{2(\pi-\alpha_{a,b})}{\pi}\displaystyle\sum_{n=21}^{\infty} \frac{1}{n} \left( \frac{2}{ \sqrt{5+\sqrt{6}}} \left( \sqrt[]{2-\sqrt[]{\frac{3}{2}}}\right) \right ) ^n \\
\end{align*}

Using the fact that, for $|a| < 1$, $\displaystyle\sum_{n=1}^\infty \frac{a^n}{n} = -\log(1-a)$, we obtain

\begin{align*}
|R_1| &\leq  \displaystyle\frac{2(\pi-\alpha_{a,b})}{\pi}\left( -\log\left( 1-\left( \frac{2}{ \sqrt{5+\sqrt{6}}} \left( \sqrt[]{2-\sqrt[]{\frac{3}{2}}}\right) \right )\right) - \displaystyle\sum_{n=1}^{20} \frac{1}{n} \left( \frac{2}{ \sqrt{5+\sqrt{6}}} \left( \sqrt[]{2-\sqrt[]{\frac{3}{2}}}\right) \right ) ^n \right) \\ 
&=  0.0000145757...  \\
&\leq 0.0000145758 .\\
\end{align*}

Therefore, since $\Omega_{a,b}(-1/2) \geq 0$, we have
 
\begin{align*}
\displaystyle \Omega_{a,b}  \left( -\frac{1}{2}\right) &= \left| \displaystyle \Omega_{a,b}\left(  -\frac{1}{2}\right) \right|  = |\Delta (a,b) + T_1 +R_1| \\
&\leq \Delta(a,b) +|T_1| +|R_1| = \Delta(a,b) + 0.3194490869562.\\
\end{align*}

Taking $t=-1/2$ in Proposition 2.4, we conclude the case $|b/a|=1$. For $|b/a|=2$, we only need to expand the series given by Proposition 2.7 until 15 terms, then, taking $t=-1$ in Proposition 2.4, we get the upper bound required. For $|b/a|=3$, the process is exactly the same but we only need to expand the series until 7 terms and take $t=-3/2$ in Proposition 2.4. Finally, the cases $|b/a| \geq 4$  and $b=0$ follow directly from Proposition 2.8. This completes the proof of the theorem. $\hfill\square$\\ \\

Theorem 2.8 gives us rigorous upper bounds for specific values of $|b/a|$. For instance, if we take $a= 7/15$, $b= 250/36$, then $b/a >4$. Therefore

\begin{center}
$\mu^{ess}(h_{a,b}) \leq \log 5 +\log 9 +\log 4 +\log (7/15) + \log(3750/252) = \log(1250).$
\end{center}

If we take the Zhang-Zagier height $h_Z = h_{1,-1}$, we have $|b/a|=1$, therefore, Theorem 2.8 gives us the upper bound
\begin{align*}
\mu^{ess}(h_Z) \leq 0.31944909.
\end{align*}

Other cases like $|b/a|= 1/2$ must be treated individually using Proposition 2.7. \\

Assume now that $a,b \in \mathbb{Q}[i] \setminus \mathbb{Q}$, we can prove the following theorem\\

\textit{{\bf{Theorem 2.10:}} Let $a,b \in \mathbb{Q}[i] \setminus  \mathbb{Q}$, then }

\begin{center}
$\mu^{ess}(h_{a,b}) \leq \Delta(a,b) + \varphi\left( \Re\left( \displaystyle \frac{b}{a}\right)\right)  +2\varphi\left(  \Im\left( \displaystyle \frac{b}{a}\right)\right) $.
\end{center}

\textit{Proof:}  Using Proposition $2.4$, we obtain that $\mu^{ess}(h_{a,b})$ is less or equal than

\begin{center}
$\Omega_{a,b}(t)=  \Delta(a,b) + \varphi(t) + \displaystyle \frac{1}{2\pi} \int_0^{2\pi} \log^{+} \left|  e^{i\theta} + \displaystyle\frac{b}{a} + t\right| + \log^{+} \left|  e^{i\theta} + \displaystyle\frac{\overline{b}}{\overline{a}} + t \right| d\theta $.
\end{center}

Since

\begin{align*}
\displaystyle \int_0^{2\pi} \log^{+} \left|  e^{i\theta} + \displaystyle\frac{\overline{b}}{\overline{a}} + t \right| d\theta  &=\displaystyle  \int_0^{2\pi} \log^{+} \left| \overline{ e^{-i\theta} + \displaystyle\frac{b}{a} + t} \right| d\theta =  \int_0^{2\pi} \log^+ \left| e^{-i\theta} +\displaystyle\frac{b}{a}+ t \right| d\theta  \\
&= \int_{-2\pi}^{0} \log^{+} \left| e^{i\theta} + \displaystyle\frac{b}{a} + t \right| d\theta = \int_0^{2\pi} \log^{+} \left|  e^{i\theta} + \displaystyle\frac{b}{a} + t\right|d\theta.\\
\end{align*}

We conclude that

\begin{center}
$\Omega_{a,b}(t)=  \Delta(a,b) + \varphi(t) + 2\left(  \displaystyle \frac{1}{2\pi} \int_0^{2\pi} \log^{+} \left|  e^{i\theta} + \displaystyle\frac{b}{a} + t\right| d\theta \right) $.
\end{center}

Now, we evaluate the function at $t= -\Re(b/a)$

\begin{center}
$\Omega_{a,b}\left( -\Re\left( \displaystyle \frac{b}{a}\right)  \right) = 2\log|a| + \varphi\left( -\Re\left( \displaystyle \frac{b}{a}\right)\right)  + 2 \left( \displaystyle \frac{1}{2\pi} \int_0^{2\pi} \log^{+} \left|  e^{i\theta} +  i \Im \left(  \displaystyle \frac{b}{a}\right) \right| d\theta\right) .$
\end{center}

Note that

\begin{align*}
\displaystyle\int_0^{2\pi} \log^{+} \left|  e^{i\theta} +  i \Im \left(  \displaystyle \frac{b}{a}\right) \right| d\theta &= \int_0^{2\pi} \log^{+} \left|  e^{i(\theta +\frac{3\pi}{2})} +  \Im \left(  \displaystyle \frac{b}{a}\right) \right| d\theta =\int_{3\pi/2}^{2\pi +3\pi/2} \log^{+} \left|  e^{i\theta} +   \Im \left(  \displaystyle \frac{b}{a}\right) \right| d\theta \\
&= \int_{3\pi/2}^{2\pi} \log^{+} \left|  e^{i\theta} +   \Im \left(  \displaystyle \frac{b}{a}\right) \right| d\theta + \int_{2\pi}^{2\pi +3\pi/2} \log^{+} \left|  e^{i\theta} +  \Im \left(  \displaystyle \frac{b}{a}\right) \right| d\theta \\
&= \int_{3\pi/2}^{2\pi} \log^{+} \left|  e^{i\theta} +   \Im \left(  \displaystyle \frac{b}{a}\right) \right| d\theta + \int_{0}^{3\pi/2} \log^{+} \left|  e^{i\theta} +  i \Im \left(  \displaystyle \frac{b}{a}\right) \right| d\theta = 2\pi \varphi \left( \Im \left(  \displaystyle \frac{b}{a}\right)\right).\\ 
\end{align*}

Finally, using Lemma $2.5$, we conclude the theorem.  $\hfill\square$ \\ \\

\textit{\textbf{Remark:} The specific values of $t$ used in the various applications of Proposition $2.4$, were suggested by numerical experiments.} \\ \\

\begin{center}
\textit{\bf{3. Lower Bounds}}
\end{center}

In this section we will compute lower bounds for $\mu^{ess}(h_{a,b})$ for $a,b \in \overline{\mathbb{Q}}$ and $a \neq 0$. We use the method described in \cite{MR3802441} section 2.2. We will assume that we are in the non-toric case. For every $\sigma \in G(a,b)$, we consider the real-valued functions $g_{\sigma},f_{\sigma},G_{\sigma}$, given by

\begin{align*}
g_{\sigma}(z) &= \log^+|z| +\log^+|\sigma(a)z+\sigma(b)|, \\
f_{\sigma}(z) &= \log^+\displaystyle \left| z \right| +\log^+\left| \frac{1}{\sigma(a)z+\sigma(b)} \right|, \\
G_{\sigma}(z) &= \log^+\displaystyle \left| z \right| +\log^+\left| \frac{\sigma(a)+\sigma(b)z}{z} \right|. \\
\end{align*}

We have that these functions go to $\infty$ when $|z| \to \infty$. Furthermore, $f_\sigma \to \infty$ when $z \to-\sigma(b)/\sigma(a)$ and $G_\sigma \to \infty$ when $z \to 0$, and they are continuous elsewhere, so they attain their minimum values. We denote by $\min(g_\sigma), \min(f_\sigma)$ and $\min(G_{\sigma})$, the minimum value of each of these functions respectively. Then, we define $g^{\min}= (1/ [K_{a,b}:\mathbb{Q}])\sum_{\sigma \in G(a,b)} \min (g_{\sigma}) $, $f^{\min}= (1/ [K_{a,b}:\mathbb{Q}])\sum_{\sigma \in G(a,b)} \min (f_{\sigma}) $ and $G^{\min}= (1/ [K_{a,b}:\mathbb{Q}])\sum_{\sigma \in G(a,b)} \min (G_{\sigma}) $. Finally, we define

\begin{center}
 $\mathcal{L}(a,b) = \max \lbrace g^{\min},f^{\min},G^{\min}\rbrace $.
\end{center}

Assume now that the minimum value of $h_{a,b}$ is achieved only at a finite non empty set of algebraic numbers. Let $\alpha_1,\alpha_2,..., \alpha_k$ be the algebraic numbers where $h_{a,b}$ is equal to the minimum. We consider $\left\lbrace  f_1,f_2,..., f_r \right\rbrace$ a set monic irreducible polynomials, such that their combined roots are $\lbrace \alpha_1,\alpha_2,... , \alpha_k \rbrace$. Now, let  $A_1,A_2,...,A_r \in \mathbb{R}_{\geq 0}$ be such that, $ A_1\deg(f_1)A_2\deg(f_2)...A_r\deg(f_r) < 2$ and for every $\alpha \in \overline{\mathbb{Q}} \setminus \bigcup_{i=1}^k \text{Gal}(\alpha_i)$, and every non-Archimedean place $\nu$ in $M_{\overline{\mathbb{Q}}}$, we have

\begin{equation}
\log^+ |\alpha|_\nu +\log^+|a\alpha+b|_\nu \geq \displaystyle \sum_{i=1}^r A_i \log |f_i(\alpha)|_{\nu}.
\end{equation}

We call $P$ the set of all $(A_1,A_2,..., A_r)\in \mathbb{R}^r$, such that these conditions hold. Since $(0,0,0,..., 0) \in P$, we have that $P \neq \emptyset$. \\

\textit{{\bf{Lemma 3.1:} } The set $P$ defined before is bounded} \\

\textit{Proof: } Firstly, we have that $A_i \geq 0$, therefore each $A_i$ is bounded from below. We fix $\nu_0 \in M_{\overline{\mathbb{Q}}}$ and consider $\alpha \in \overline{\mathbb{Q}} \setminus \bigcup_{i=1}^k \text{Gal}(\alpha_i)$, such that, $\log^+|a\alpha+b|_{\nu_0} = \log|a\alpha|_{\nu_0}$ and for each $1 \leq i \leq r$,  $|f_i(\alpha)|_{\nu_0} = |\alpha^{\deg(f_i)}|_{\nu_0} \geq 1$. Therefore, from (8) we have   $\log|a\alpha^2|_{\nu_0} \geq  \sum_{i=1}^r A_i\log|\alpha^{\deg(f_i)}|_{\nu_0} \geq A_{t}\log|\alpha^{\deg(f_t)}|_{\nu_0}= \log|\alpha^{\deg(f_t)A_t}|_{\nu_0}$, for each $t \in \lbrace 1,2,... ,r\rbrace $. For $|\alpha|_{\nu_0}$ large enough, this equality holds only if $A_t \leq 2/\deg(f_t)$. This completes the proof of the lemma. $\hfill\square$ \\\\

Now, for each $\sigma \in G(a,b)$ we define the real valued function $g_{A_1,...,A_r,\sigma}$ by

\begin{center}
$g_{A_1,...,A_r,\sigma}(z) =  \log^+ |z| +\log^+|\sigma(a)z+\sigma(b)| - \displaystyle \sum_{i=1}^r A_i \log |f_i(z)|$.
\end{center}

Since $A_i \geq 0$ for each $i \in \lbrace 1,2,..., r \rbrace$, and $A_1\deg(f_1)A_2\deg(f_2)...A_r\deg(f_r) < 2$, we have that $ g_{A_1,...,A_k,\sigma} \to \infty$ when $|z| \to \infty$ or $z \to x$, where $x \in \bigcup_{i=1}^k \text{Gal}(\alpha_i)$, and it is continuous elsewhere, so it attains its minimum value. Now, we consider the function $H_{\sigma}: P \to \mathbb{R}$, given by

\begin{center}
$H_{\sigma}(A_1,A_2,..., A_r) = \displaystyle \inf g_{A_1,...,A_r,\sigma} $.
\end{center}

We define 

\begin{center}
$\tau(a,b)=\displaystyle \frac{1}{[K_{a,b}:\mathbb{Q}]} \displaystyle\sup_{(A_1,...,A_r) \in P}\sum_{\sigma \in G(a,b)} H_{\sigma}(A_1,...,A_r)$. 
\end{center}

Now, we can state the following theorem. \\

\textit{{\bf{Theorem 3.2:}} Let $a,b \in \overline{\mathbb{Q}}$ with $a \neq 0$. Then, }

\begin{center}
$\mu^{ess} (h_{a,b}) \geq \mathcal{L}(a,b)$.
\end{center}

\textit{Moreover, if the minimum value is achieved only at a finite non empty set of algebraic numbers, then}
\begin{center}
$\mu^{ess} (h_{a,b}) \geq \tau(a,b)$.
\end{center}

Before proving the theorem we need the following lemma \\

\textit{{\bf{Lemma 3.3:}} Let $f: \overline{\mathbb{Q}} \to \mathbb{R}$. Define $f^{in} :\overline{\mathbb{Q}}^* \to \mathbb{R}$ by $f^{in} (\alpha) =f(1/\alpha)$. Then }

\begin{center}
$\mu^{ess} (f) =\mu^{ess}(f^{in})$.
\end{center}

\textit{Proof:} We write, $\mu^{ess}(f) = M$. Then, there exists a sequence of distinct algebraic numbers $\lbrace \gamma_n \rbrace_{n \in \mathbb{N}}$ such that $f(\gamma_n) \to M$. In particular, we may assume $\gamma_n \neq 0$ for all $n \in \mathbb{N}$. Then $\lbrace \beta_n \rbrace_{n \in \mathbb{N}}$ given by $\beta_n= 1/\gamma_n$ is a sequence of algebraic numbers, and  $f^{in}(\beta_n)=f(\gamma_n) \to M$, therefore $\mu^{ess} ({f^{in}}) \leq M$. The other  inequality is similar. $\hfill\square$\\

\textit{Proof of Theorem 3.2: } Let $K /\mathbb{Q}$ be a Galois extension such that $a,b,\alpha \in K$. Note that

\begin{align*}
h_{a,b}(\alpha) &=\displaystyle \frac{1}{[K :\mathbb{Q}]} \left( \sum_{p}\sum_{ \nu | \left| . \right|_p } \log^+|\alpha|_\nu+\log^+ \left| a\alpha +b \right|_\nu  + \sum_{\nu | |.|} \log^+|\alpha|_\nu +\log^+ \left| a\alpha +b \right|_\nu \right)\\
&\geq \displaystyle \frac{1}{[K :\mathbb{Q}]} \sum_{\nu | |.|} \log^+|\alpha|_\nu +\log^+ \left| a\alpha +b \right|_\nu  \\
&=  \displaystyle \frac{1}{[K :\mathbb{Q}]} \sum_{\sigma \in G(a,b)} \sum_{\substack{ \tau \in \text{Gal}(K\setminus\mathbb{Q}) \\ \tau | \sigma}} \log^+|\tau(\alpha)| +\log^+ \left| \sigma(a)\tau(\alpha) +\sigma(b) \right| \\
&\geq \frac{1}{[K_{a,b}:\mathbb{Q}]}\displaystyle \sum_{\sigma \in G(a,b)}\inf  g_\sigma\\
&= g^{\min}.
\end{align*}

Therefore, $g^{\min} \leq h_{a,b}(\alpha)$ for all $\alpha$, we conclude that $g^{\min} \leq \mu^{ess}(h_{a,b})$. Similarly, we can conclude that $f^{\min} \leq \mu^{ess}(h_{a,b})$. In fact, we can consider the real-valued function $j_{a,b}$, given by $j_{a,b}(\alpha)=h(\alpha) + h( 1/(a\alpha+b))$. Since $h(\alpha)=h(1/ \alpha)$, we have $j_{a,b}=h_{a,b}$. However, the Archimedean parts of them are different, so, we can use the same inequalities as before and conclude that $f^{\min} \leq \mu^{ess}(j_{a,b}) = \mu^{ess}(h_{a,b}) $. Finally, let the real valued function $n_{a,b}$, be given by $n_{a,b}(\alpha)= h(\alpha) + h((a+b\alpha)/\alpha )= h_{a,b}^{in}$. Using Lemma 3.3 and the same method used before, we conclude  that $G^{\min} \leq \mu^{ess}(n_{a,b}) = \mu^{ess}(h_{a,b}^{in}) = \mu^{ess}(h_{a,b}) $. Therefore,  $\mathcal{L}(a,b) \leq \mu^{ess}(h_{a,b})$.\\

Assume now that the minimum value of $h_{a,b}$  is achieved only at a finite set $\left\lbrace \alpha_1, \alpha_2,.. , \alpha_k \right\rbrace $. Let $ \alpha \in \overline{\mathbb{Q}}\setminus \bigcup_{i=1}^k \text{Gal}(\alpha_i)$ and $A_1,A_2,...,A_r \in P$. The product formula gives us that for each $i \in \lbrace 1,2,...,r\rbrace$, we have

\begin{center}
$\displaystyle\sum_{\nu \in M_{K}} A_i \log |f_i(\alpha)|_\nu = 0$
\end{center}

Therefore

\begin{align*}
h_{a,b}(\alpha) &=\displaystyle \frac{1}{[K:\mathbb{Q}]}  \left( \sum_{p}\sum_{ \nu | \left| . \right|_p } \log^+|\alpha|_\nu+\log^+ \left| a\alpha +b \right|_\nu   - \sum_{i=1}^r A_i\log|f_i(\alpha)|_\nu\right) \\
& + \displaystyle \frac{1}{[K:\mathbb{Q}]}  \left( \sum_{\nu | |.|} \log^+|\alpha|_\nu + \log^+|a\alpha+b|_\nu -\sum_{i=1}^rA_i \log |f_i(\alpha)|_\nu\right) \\
&\geq \displaystyle \frac{1}{[K:\mathbb{Q}]}  \left( \sum_{\nu | |.|} \log^+|\alpha|_\nu + \log^+|a\alpha+b|_\nu -\sum_{i=1}^rA_i \log |f_i(\alpha)|_\nu\right)  \\
&= \displaystyle \frac{1}{[K :\mathbb{Q}]} \sum_{\sigma \in G(a,b)} \sum_{\substack{ \tau \in \text{Gal}(K/\mathbb{Q}) \\ \tau | \sigma}} g_{A_1,A_2,..., A_r, \sigma} (\tau(\alpha)) \\
&\geq \frac{1}{[K_{a,b}:\mathbb{Q}]}\displaystyle \sum_{\sigma \in G(a,b)}H_\sigma(A_1,A_2,... , A_r).\\
\end{align*}

Since the last inequality holds for all $(A_1,A_2,...,A_r) \in P$, we have that

\begin{align*}
h_{a,b}(\alpha)&\geq \frac{1}{[K_{a,b}:\mathbb{Q}]}\sup_{(A_1,A_2,..., A_k) \in P}\displaystyle \sum_{\sigma \in G(a,b)} H_\sigma(A_1,A_2,... , A_r)\\
&= \tau(a,b).
\end{align*}

Since the last inequality holds  for all $\alpha$ except finitely many, we conclude that $\tau(a,b) \leq \mu^{ess}(h_{a,b})$. $\hfill\square$ \\ \\

The following observation will be useful to compute lower bounds. Since $g_{A_1,...,A_k,\sigma}$ is harmonic off the two sets $|z|=1$ and $|\sigma(a)z+\sigma(b)|=1$, then the minimum is achieved only on these sets. The same happens with  $f_{\sigma}$ and $g_{\sigma}$. For $G_{\sigma}$, the function is harmonic off the two sets $|z|=1$ and $|(\sigma(a)+\sigma(b)z)/z|=1$, therefore the minimum is achieved on these sets. Now we are ready to prove the following proposition \\ \\

\textit{{\bf{Proposition 3.4:}} Let $a,b \in \overline{\mathbb{Q}}$ with $a \neq 0$. Assume that, there exists $\sigma_0 \in G(a,b)$, such that $|\sigma_0(b)|-|\sigma_0(a)| > 1$. Then}
\begin{center}
$\displaystyle \frac{1}{[K_{a,b}: \mathbb{Q}]}\min \left( \displaystyle \log(|\sigma_0(b)|-|\sigma_0(a)|), \log \left(\frac{|\sigma_0(b)|-1}{|\sigma_0(a)|}\right)\right) \leq \mu^{ess} (h_{a,b})$.
\end{center}

\textit{Proof: } Consider the function

\begin{center}
$g_{\sigma_0}(z) = \log^+|z|+\log^+|\sigma_0(a)z+\sigma_0(b)|$.
\end{center}

Firstly, we consider $z=e^{i\theta}$, with $0\leq \theta \leq 2\pi$, $\arg(\sigma_0(a))=\gamma$ and $\arg(\sigma_0(b))=\beta$, note that $|\sigma_0(a)e^{i\theta} + \sigma_0(b)| > 1$, so we have

\begin{align*}
g_{\sigma_0}(e^{i\theta}) &= \log|\sigma_0(a)e^{i\theta}+\sigma_0(b)| \\
&= \displaystyle \frac{1}{2}\log (|\sigma_0(a)|^2+|\sigma_0(b)|^2 + 2|\sigma_0(b)||\sigma_0(a)|\cos(\theta+\gamma-\beta)).
\end{align*}

The minimum value is achieved when $\cos(\theta+\gamma-\beta)= -1$, and the minimum value is $\log(|\sigma_0(b)|-|\sigma_0(a)|)$. On the other hand, note that $(|\sigma_0(b)|-1)/|\sigma_0(a)| > 1$, therefore, $|e^{i\theta}-\sigma_0(b)|/|\sigma_0(a)| > 1$. We conclude that

\begin{align*}
g_{\sigma_0}\left( \displaystyle \frac{e^{i\theta}-\sigma_0(b)}{\sigma_0(a)} \right) &= \log \left| \frac{e^{i\theta}-\sigma_0(b)}{\sigma_0(a)} \right| \\
&= \displaystyle \frac{1}{2}\log \left( \frac{|\sigma_0(b)|^2+1-2|\sigma_0(b)|\cos(\theta-\beta)}{|\sigma_0(a)|^2} \right).
\end{align*}

The minimum is achieved when $\theta=\beta$, and its value is $\log ((|\sigma_0(b)|-1|)/|\sigma_0(a)|)$. Since, $\mathcal{L}(a,b) \geq \min g_{\sigma_0}/[K_{a,b}:\mathbb{Q}]$, using Theorem 3.2  we conclude the proof of the proposition. $\hfill\square$ \\ \\

\textit{{\bf{Proposition 3.5:}} Let $a,b \in \overline{\mathbb{Q}}$, $b \neq 0$. Assume that there exist $\sigma_0 \in G(a,b)$, such that, $|\sigma_0(a)|-|\sigma_0(b)|| > 1$. Then}
\begin{center}
$\displaystyle \frac{1}{[K_{a,b}:\mathbb{Q}]}\log \left(\displaystyle\frac{|\sigma_0(a)|}{|\sigma_0(b)|+1}\right) \leq \mu^{ess}(h_{a,b})$.
\end{center}

\textit{Proof:} we consider 

\begin{center}
$G_{\sigma_0}(z)= \log^+|z|+\log^+\left| \displaystyle \frac{\sigma_0(a) + \sigma_0(b)z}{z}\right|$.
\end{center}

Let's consider, $z=e^{i\theta}$, and let $\arg(\sigma_0(a))= \gamma$ and $\arg(\sigma_0(b))= \beta$. Since $|\sigma_0(a)e^{i\theta}+\sigma_0(b)| >1$,

\begin{align*}
G_{\sigma_0}(e^{i\theta}) &= \log| \sigma_0(b)e^{i\theta} + \sigma_0(a)| \\
&= \frac{1}{2} \log (|\sigma_0(a)|^2+|\sigma_0(b)|^2 + 2|\sigma_0(b)||\sigma_0(a)|\cos(\theta+\beta-\gamma)).
\end{align*}

The minimum value is achieved when $\cos(\theta+\beta-\gamma)= -1$, an the value is $\log(|\sigma_0(a)|-|\sigma_0(b)|)$. On the other hand, note that $|\sigma_0(a)|/(|\sigma_0(b)|+1) > 1$, therefore, $|\sigma_0(a)|/|e^{i\theta}-\sigma_0(b)| > 1$. We conclude that

\begin{align*}
G_{\sigma_0}\left( \displaystyle \frac{\sigma_0(a)}{e^{i\theta}-\sigma_0(b)} \right) &= \log \left| \frac{\sigma_0(a)}{e^{i\theta}-\sigma_0(b)} \right| \\
&= -\displaystyle \frac{1}{2}\log \left( \frac{|\sigma_0(b)|^2+1-2|\sigma_0(b)|\cos(\theta-\beta)}{|\sigma_0(a)|^2} \right).
\end{align*}

The minimum is achieved when $\cos(\theta-\beta)= -1$, and its value is $\log (|\sigma_0(a)|/(|\sigma_0(b)|+1))$. Note that, since $|\sigma_0(a)|-|\sigma_0(b)| > 1$, it is not hard to prove that $|\sigma_0(a)|-|\sigma_0(b)| > |\sigma_0(a)|/(|\sigma_0(b)|+1)$. Since, $\mathcal{L}(a,b) \geq G^{\min} \geq [K_{a,b}/\mathbb{Q}]\min G_{\sigma_0}$, using Theorem 3.2  we conclude the proof of the proposition. $\hfill\square$ \\

\textit{{\bf{Proposition 3.6:}} Let $a,b \in \overline{\mathbb{Q}}$, $a \neq 0$, $b \neq 0$. Assume that, there exist $\sigma_0 \in G(a,b)$, such that, $0 \leq | |\sigma_0(a)|-|\sigma_0(b)|| \leq 1$ and $|\sigma_0(a)|+|\sigma_0(b)| < 1$. Then}
\begin{center}
$\displaystyle \frac{1}{[K_{a,b}: \mathbb{Q}]}\log \left( \displaystyle\frac{1}{|\sigma_0(a)|+|\sigma_0(b)|} \right) \leq \mu^{ess}(h_{a,b})$.
\end{center}

\textit{Proof: } In this case, we consider

\begin{center}
$f_{\sigma_0}(z) = \log^+\left| \displaystyle z\right| +\log^+\left| \displaystyle \frac{1}{\sigma_0(a)z +\sigma_0(b)}\right|$.
\end{center}

If $z=e^{i\theta}$, then $1/|\sigma_0(a)z+\sigma_0(b)| >1$, therefore

\begin{center}
$f_{\sigma_0}(e^{i\theta}) = \log\left| \displaystyle \frac{1}{\sigma_0(a)e^{i\theta} +\sigma_0(b)}\right|$.
\end{center}

The minimum value is $ \log(1/||\sigma_0(a)|+|\sigma_0(b)||)$, now, assume $z= (1-\sigma_0(b)e^{i\theta})/(\sigma_0(a)e^{i\theta})$, we note that $|(1-\sigma_0(b)e^{i\theta})/(\sigma_0(a)e^{i\theta})| \geq (1-|\sigma_0(b)|)/|\sigma_0(a)| > 1$, therefore

\begin{center}
$f_{\sigma_0} \left( \displaystyle \frac{1-\sigma_0(b)e^{i\theta}}{\sigma_0(a)e^{i\theta}} \right) = \log\left| \displaystyle \frac{1-\sigma_0(b)e^{i\theta}}{\sigma_0(a)e^{i\theta}}\right|$.
\end{center}

The minimum value is $\log((1-|\sigma_0(b)|)/|\sigma_0(a)|)$. Since $|\sigma_0(a)|+|\sigma_0(b)| < 1$, it is not hard to prove that  $(1-|\sigma_0(b)|)/|\sigma_0(a)| > 1/(|\sigma_0(a)|+|\sigma_0(b)|)$ , we conclude that the minimum value of $f_{\sigma_0}$ is $\log (1/(|\sigma_0(a)|+|\sigma_0(b)|))$. Since, $\mathcal{L}(a,b) \geq f^{\min} \geq [K_{a,b}/\mathbb{Q}] \min f_{\sigma_0}$, using Theorem 3.2  we conclude the proof of the proposition.. $\hfill\square$ \\

Propositions $3.4$, $3.5$ and $3.6$ give us non zero lower bounds for  the cases $a,b \in  \overline{\mathbb{Q}}$, such that there exists $\sigma_0 \in G(a,b)$, such that $||\sigma_0(a)|-|\sigma_0(b)|| > 1$, $b \neq 0$, or the cases where $||\sigma_0(a)|-|\sigma_0(b)|| \leq 1$ and $|\sigma_0(a)|+|\sigma_0(b)| < 1$. Using Theorem 3.2 together with Propositions 2.1, 3.4, 3.5 and 3.6  we obtain Theorem B (Note that, in Theorem B, we have used $\mathcal{K}(a,b)$, this number can represent $\mathcal{L}(a,b)$ or $\tau(a,b)$).\\

\textit{{\bf{Lemma 3.7:}} Let $a,b \in \overline{\mathbb{Q}}$, the following sentences are equivalent } \\
\begin{enumerate}[(a)]
\item[i)] \textit{The minimum value of $h_{a,b}$ is achieved and $\min(h_{a,b}) = 0$.}
\item[ii)] \textit{$b=0$ or \textit{$b$ is a root of unity, or there exists a root of unity $\zeta$, such that $a\zeta +b=0$ or $a\zeta +b=\zeta_0$ a root of unity}.}
\end{enumerate}

\textit{Furthermore, assume that any of this sentences hold, if $b \neq 0$, $b$ is not a root of unity and $|a| \neq |b|$, then necessarily $||a|-|b|| \leq 1$ and $|a| + |b| \geq 1$}\\

\textit{Proof: } (i)$ \Rightarrow$(ii) Assume that $\min h_{a,b}$ is achieved at $\alpha \in \overline{\mathbb{Q}}$. Since, $h_{a,b}(\alpha)=h(\alpha)+h(a\alpha+b)$, we need $h(\alpha)=0$ and $h(a\alpha +b)=0$. The first equality implies $\alpha=0$ or $\alpha= \zeta$ a root of unity. If $\alpha=0$, we also need $h(b)=0$, so there are two options, $b=0$ or $b$ is a root of unity. On the other hand, if $\alpha=\zeta$ is a root of unity, we also need $a\zeta+b=0$ or $a\zeta+b$  a root of unity.\\

(ii)$ \Rightarrow$(i) If $b=0$ or a root of unity, we take $\alpha=0$ and we obtain $h_{a,0}(\alpha)=0$. On the other hand, if there exists a root of unity $\zeta$ such that  $a\zeta +b=0$ or $a\zeta +b=\zeta_0$ a root of unity, we take $\alpha = \zeta$, and we obtain $h_{a,b}(\alpha)=0$. $\hfill\square$  \\

Now, assume that $b \neq 0$, $b$ is not a root of unity and $|a| \neq |b|$, by (ii), we have that there exist two roots of unity $\zeta$ and $\zeta_0$, such that $a\zeta + b = \zeta_0$, therefore $|a\zeta +b| = 1$, we conclude that $||a|-|b|| \leq |a\zeta+b| = 1$ and $1 =|a\zeta +b| \leq |a\zeta| +|b| = |a| +|b|$. $\hfill\square$ \\ \\

\textit{{\bf{Corollary 3.8:}} Let $a \in \overline{\mathbb{Q}}\setminus  \lbrace 0 \rbrace$. Then $\mu^{ess}(h_{a,0})=0$ if and only if $a$ is a root of unity. Moreover, }

\begin{center}
$\mu^{ess}(h_{a,0}) \geq h(a)$
\end{center}

\textit{Proof: }  We consider $a \in \overline{\mathbb{Q}} \setminus \lbrace 0 \rbrace$ not a root of unity, then $h_{a,0}(0)=0$, and it is zero only at $0$, therefore $f_1(x)=x$.  Firstly we need to determine the possible values of $A_1$. For each $\nu$ non-Archimedean and $\alpha \in \overline{\mathbb{Q}}\setminus \lbrace 0 \rbrace$, we need that

\begin{center}
$\log^+|\alpha|_\nu +\log^+|a\alpha|_\nu \geq A_1\log|\alpha|_\nu$.
\end{center}

if $|\alpha|_\nu \leq 1$, this inequality always occur. If $|\alpha|_\nu > 1$ and $|a\alpha|_\nu \leq 1$,  we get the restriction $A_1 \leq 1$. On the other hand, if $|\alpha|_\nu < 1$ and $|a\alpha|_\nu < 1$, we get the already known restriction $0 \leq A_1$, therefore $0 \leq A_1 \leq 1$.\\

Now, we will consider $\sigma \in G(a)$. Then,

\begin{center}
$g_{A_1, \sigma}(z)= \log^+|z|+\log^+\left| \displaystyle\sigma(a)z\right| - A_1\log|z|$.
\end{center}

if $z=e^{i\theta}$, we have

\begin{center}
$g_{ A_1, \sigma}(e^{i\theta})= \log^+|\sigma(a)|$.
\end{center}

On the other hand, if $z= e^{i\theta}/\sigma(a)$, then

\begin{center}
$g_{ A_1, \sigma} \left( \displaystyle \frac{e^{i\theta}}{\sigma(a)} \right)= \log^+\displaystyle \frac{1}{|\sigma(a)|}-A_1\log\displaystyle \frac{1}{|\sigma(a)|}    $.
\end{center}

Assume that $|\sigma(a)| \geq 1$. Then, the minimum value is $H_\sigma(A_1) = A_1\log|\sigma(a)|$. On the other hand, if $|\sigma(a)| < 1$, then the minimum value is $0$, therefore  $\inf g_{A_1,\sigma} = 0$. Summarizing,  $H_\sigma(A_1) = A_1\log^+|\sigma(a)|$, therefore 

\begin{align*}
\tau(a,0) &= \displaystyle \frac{1}{[K_a:\mathbb{Q}]} \sup_{0 \leq A_1 \leq 1}\sum_{\sigma \in G(a)} A_1\log^+|\sigma(a)| \\
&= \displaystyle \frac{1}{[K_a:\mathbb{Q}]} \sum_{\sigma \in G(a)} \log^+|\sigma(a)| 
\end{align*}

Furthermore, we can improve this result, note that if we go back to the proof of theorem $3.3$, we can get that, a lower bound for $\mu^{ess}(h_a)$ is given by

\begin{center}
$ \displaystyle     \frac{1}{[K_a:\mathbb{Q}]}  \sum_{\substack{p \\ \text{prime}}}\sum_{\sigma \in G(a)} \min_{\alpha \in \overline{\mathbb{Q}}} \left( \log^+ |\alpha|_p + \log^+|\sigma(a) \alpha)|_p- \log|\alpha|_p \right) + \tau(a,0)$
\end{center}

Now, for each $p$ prime and $\sigma \in G(a)$, we set the function $f_p^{\sigma}: \overline{\mathbb{Q}} \to \mathbb{R}_{\geq 0}$, defined by $f_p^{\sigma}(\alpha)=\log^+ |\alpha|_p + \log^+|\sigma(a)\alpha|_p-\log|\alpha|_p$, we will find the minimum value of $f_p^{\sigma}$. Firstly, suppose that $|\sigma(a)|_p \leq 1$ then, we can take $\alpha$ a root of unity, and we get, $f_p^{\sigma}(\alpha)= 0$, since $f_p^\sigma \geq 0$, this will be the minimum value. On the other hand, if $|\sigma(a)|_p > 1$, let $\alpha \in \overline{\mathbb{Q}}$, if $|\alpha|_p \geq 1$, then $f_p^{\sigma}(\alpha)=\log|\sigma(a)\alpha|_p \geq \log|\sigma(a)|_p$. Assume now that, $|\alpha|_p < 1$, then $f_p^{\sigma}(\alpha)=\log^+|\sigma(a)\alpha|_p - \log|\alpha|_p$, if $|\sigma(a)\alpha|_p \leq 1$, then $f_p^{\sigma}(\alpha)= \log 1/|\alpha|_p \geq \log|a|_p$. On the other hand, if $|\sigma(a)\alpha|_p > 1$, then $f_p^{\sigma}(\alpha)=\log|\sigma(a)\alpha|_p - \log|\alpha|_p=\log|\sigma(a)|_p$. Summarizing, we have that $\min (f_p^\sigma) = \log^+ |\sigma(a)|_p$ . We conclude that 

\begin{center}
$\mu^{ess} (h_{a,0}) \geq \displaystyle \frac{1}{[K_a:\mathbb{Q}]}  \sum_{\substack{p \\ \text{prime}}}\sum_{\sigma \in G(a)} \log^+ |\sigma(a)|_p + \tau(a,0)=h(a)$
\end{center}

This concludes the proof of the corollary. $\hfill\square$ \\\\

Note that, using $(i)$ of Proposition $2.7$, we get that $\mu^{ess} (h_{a,0}) \leq h(a)$, therefore, we have that $\mu^{ess}(h_{a,0})= h(a)$. \\ \\

\textit{{\bf{Corollary 3.9:}}} $\mu^{ess}(h_{1,2}) \geq \log(\sqrt{3})$. \\

\textit{Proof: } Note that, $\alpha$ and $\alpha + 2$ are both roots of unity only at $\alpha=-1$, therefore, $f_1(x)= x+1$.  Let $A_0,A_1 \in \mathbb{R}$, for each $\alpha \in \overline{\mathbb{Q}}$ and $p$ prime, we need 

\begin{center}
$\log^+|\alpha|_p + \log^+|\alpha+2|_p \geq A_1 \log|\alpha +1 |_p$.
\end{center}

If $|\alpha|_p \leq 1$, then we need $A_1 \geq 0$. On the other hand, if $|\alpha|_p > 1$, we need $A_1 \leq 2$, so we need $0 \leq A_1 <2$. Now, we take the function

\begin{center}
$g_{A_1}(z)=\log^+|z|+\log^+|z+2|-A_1\log|z+1|.$
\end{center}

The minimum is achieved when $z=e^{i\theta}$ or $z=e^{i\theta}-2$. If $z=e^{i\theta}$, we have

\begin{align*}
g_{A_1}(e^{i\theta}) &= \log|e^{i\theta}+2| -A_1\log|e^{i\theta}+1| \\
&= \displaystyle  \frac{1}{2} \log (5+4\cos(\theta)) - \frac{A_1}{2} \log (2+2\cos(\theta)) .
\end{align*}

Taking the derivative and equalizing to zero, we get the following possible values of theta

\begin{center}
$\theta = 0$, or $\cos(\theta) = \displaystyle \frac{A_1}{4(1-A_1)} - 1 $.
\end{center}
 
On the other hand, if $z=e^{i\theta} -2$, then

\begin{align*}
g_{A_1}(e^{i\theta}-2) &= \log|e^{i\theta} -2| - A_1\log|e^{i\theta} -1 |   \\
&= \displaystyle   \frac{1}{2} \log (5-4\cos(\theta))  - \frac{A_1}{2} \log (2-2\cos(\theta)).
\end{align*}

Again, taking the derivative and equalizing to zero, we get the next possible values of theta

\begin{center}
$\theta = \pi$, or $\cos(\theta) = \displaystyle 1 - \frac{A_1}{4(1-A_1)}$.
\end{center}

It is not hard to see that the minimum value of the function will be the same in both sets, so we can take anyone. In order to find the minimum vale, we have to know if the function is lower at $\theta= 0$ or at $\cos(\theta)=A_1/(4(1-A_1))-1$. Firstly, we will assume $0 \leq A_1 <1$, evaluating at these two values we obtain $\log(9)-A_1\log(4)$ and $ \log(1/(1-A_1))- A_1\log(A_1/(2(1-A_1)))$, we can take these two values and consider the function $l:[0,1) \to \mathbb{R}$, given by $l(x)= \log(1/(1-x))- x\log(x/(2(1-x))) -(\log(9)-x\log(4)) $. Taking the derivative of $l$, we obtain $l'(x)= -\log(x/(8(1-x)))$ and $l''(x)= -1/(x(1-x))$. Therefore, we conclude that, $x=8/9$ is an absolute maximum of $l$. Since $l(8/9) =0$, we conclude that the minimum value is $\log(1/(1-A_1)) - A_1\log(A_1/(2(1-A_1))) $, and it is achieved at  $\cos(\theta)=A_1/(1-A_1)-1$. Now, we consider the function respect to $A_1$, $H_{id}(A_1)= \log(1/(1-A_1)) - A_1\log(A_1/(2(1-A_1)))$, taking the derivative respect to $A_1$ and equalizing to zero, we obtain $\log(A_1/2(1-A_1))=0$, therefore $A_1/2(1-A_1)=1$, and we conclude that the maximum value of $H_{id}(A_1)$ is achieved at $A_1=2/3$, replacing this, we obtain $H_{id}(2/3) = \log(\sqrt{3})$. On the other hand, if $1 \leq A_1 <2$, then, there is no $\theta$ such that $\cos(\theta) = A_1/(4(1-A_1)) -1$. We conclude that, the minimum value of $g_{A_1}$ is $(1/2)(\log(9)-A_1\log(4))$, the maximum value is achieved at $A_1 = 1$ and $H_{id}(1) = \log (3/2)$. Since $ \sqrt{3} >3/2$, we have $\tau(a,b)= \displaystyle \sup_{0 \leq A_1 <2} H_{id} (A_1) = \sqrt{3}$, therefore $\log(\sqrt{3}) \leq \mu^{ess}(h_{1,2})$. \\ \\

\begin{center}
\textit{\bf{4. Intervals of density}}
\end{center}

The goal of this section is to determine intervals where the image of $h_{a,b}$ is dense. Let $v \in M_\mathbb{Q}$, we denote by  $\mathbb{Q}_v$ the completion of $\mathbb{Q}$ at $v$. Let $\overline{\mathbb{Q}}_{v}$, the algebraic closure of $\mathbb{Q}_v$, and let $\mathbb{C}_{v}$ denote the completion of $\overline{\mathbb{Q}}_{v}$. It is well known that $\mathbb{C}_{v}$ is algebraically closed. Using this notation, a Galois invariant adelic set is a set of the form \\

\begin{center}
$\mathbb{E}= \displaystyle \prod_{v \in \mathbb{Q}} E_v$.
\end{center}

Where $E_v$ is a subset of of $\mathbb{C}_{v}$ invariant under the action of the absolute $v$-adic Galois group Gal$(\overline{\mathbb{Q}}_{v}/\mathbb{Q})$, for all $v$, and such that  $E_v = \mathcal{O}_v$ for all but a finite number of $v$. The capacity of $\mathbb{E}$ is defined by $\text{Cap}(\mathbb{E})= \prod_{v} \text{Cap}(E_v)$, it is well defined, because $\text{Cap}(E_v)= \text{Cap}(\mathcal{O}_v)=1$ (see \cite{MR1009368}, section $4$), for all but a finite number of $v$. Furthermore, for each non-Archimedean $v$, write $B_v(0,r)= \left\lbrace z \in \mathbb{C}_v : |z|_v < r \right\rbrace$, then, the equilibrium measure of $\overline{B}_v(0,r)$ is the linear functional $\mu_{\overline{B}_v(0,r)}: C^0(\overline{B}_v(0,r),\mathbb{R}) \to \mathbb{R}$, defined by

\begin{center}
$\displaystyle\int_{\overline{B}_v(0,r)} f(t)d\mu_{\overline{B}_v(0,r)}(t) = \sup_{|z| \leq r}|f(z)|$
\end{center}

(see \cite{MR1009368}, section $4$, for more details). The following proposition will be useful \\

\textit{{\bf{Proposition 4.1: }}Let $\mathbb{E}= \prod_{v \in M_\mathbb{Q}} E_v$ an adelic set with $\text{Cap}(\mathbb{E})=1$. Then, there exists a sequence $(x_l)_{l \in \mathbb{N}}$ of pairwise distinc points of $\overline{\mathbb{Q}}^{\times}$, with $Gal(x_l)_v \subset B(E_v, 1/l)$, for all $v \in M_\mathbb{Q}$. Furthermore, for all $v \in M_\mathbb{Q}$, the sequence of measures $(\delta(\text{Gal}(x_l,v)))_{l \in \mathbb{N}}$ converges to the equilibrium measure $\mu_{\overline{E}_v}$ weakly. } \\

\textit{Proof: } Direct from \cite{MR3928039} Proposition $7.4$.\\

 Let $a,b \in \mathbb{Q}$, we consider $S_{a,b}= \lbrace p_1, p_2, ..., p_s\rbrace$ and $\mathbb{Q}^+ = \lbrace q \in \mathbb{Q} : q > 0 \rbrace$, then, we define  $\Gamma_{a,b}: \mathbb{R}\times (\mathbb{Q}^+)^s \to \mathbb{R}$, given  by

\begin{align*}
\Gamma_{a,b} (x,r_1,r_2,...,r_s)= \displaystyle \sum_{i=1}^s \log^+ r_i+\log^+ \max \lbrace|a|_{p_i} r_i, |b|_{p_i}\rbrace + \frac{1}{2\pi} \int_0^{2\pi} \log^+ \left| \frac{e^{i\theta}}{r_1r_2... r_s} +x \right| + \log^+\left|\frac{ae^{i\theta}}{r_1r_2... r_s}+b+ax \right| d\theta.
\end{align*}

\textit{{\bf{Theorem 4.2:}} Let $a,b \in \mathbb{Q}$, then for each $x \in \mathbb{R}$ and $r_1,r_2,..., r_s \in \mathbb{Q}^+$, the image of $h_{a,b}$ is dense in the interval $[\Gamma_{a,b}(x,r_1,r_2,..., r_s), \infty)$.  Furthermore, 	}
\begin{align*}
\mu^{ess} (h_{a,b}) \leq \Gamma_{a,b}(x,r_1,r_2,..., r_s).
\end{align*}
\hspace{1cm}

\textit{Proof: } 
We use Proposition $4.1$, for  $p_i \in S_{a,b}$ we take $E_{p_i}= B_{p_i}(0,r_i)= \left\lbrace z \in \mathbb{C}_{p_i} : |z|_{p_i} < r_i \right\rbrace $, where $r_i \in \mathbb{Q}^+$. For $p \notin S_{a,b}$ we take, $E_{p}= \mathcal{O}_{p}$. For $|.|_\nu= |.|_\infty$, we take $E_{|.|}= B(x, 1/(r_1r_2 ... r_s))$, where $x \in \mathbb{R}$. We define, $\mathbb{E}= \prod_{v}E_v$, then $\mathbb{E}$ is an adelic set such that $\text{Cap}(\mathbb{E})=1$.  Therefore, by Proposition $4.1$, there is a sequence $(\alpha_n)_{n \in \mathbb{N}}$ satisfying  $\text{Gal}(\alpha_{n,\nu}) \subset B (E_\nu, 1/n)$, for each $\nu \in M_{\mathbb{Q}}$, such that $\delta(\text{Gal}(\alpha_{n},\nu)) \overset{*}{\longrightarrow} \mu_{\overline{E}_\nu}$, therefore

\begin{equation}
h_{a,b}(\alpha_n) \underset{n \to \infty}{\longrightarrow} \displaystyle \sum_{\nu \in M_\mathbb{\mathbb{Q}} } \int_{\overline{E}_\nu} \log^+ |t|_\nu + \log^+|at+b|_\nu d\mu_{\overline{E}_\nu}(t) := M_{E_x}.
\end{equation}

We claim that $ M_{E_x} < \infty$. Let $p$ a prime number and $d \in \mathbb{R}$, then

\begin{center}
$\displaystyle \sup_{|z|_p \leq d}  \log^+ |az+b|_p \leq \sup_{|z|_p \leq d}  \log^+ \max \lbrace |a|_p|z|_p,|b|_p\rbrace = \log^+ \max \lbrace |a|_pd, |b|_p\rbrace$
\end{center}

If $|a|_pd \leq |b|_p$, then, taking $t_0$, such that $|t_0|_p < d$, we obtain that $|at_0+b|_p=|b|_p= \max\lbrace |a|_pd,|b|_p \rbrace$. If $|a|_pd > |b|_p$, we can take $t_1$, such that $|t_1|_p =d$, we obtain that $|at_1+b|_p=|a|_pd$. Therefore $\sup_{|z|_p \leq d}  \log^+ |az+b|_p  = \log^+ \max \lbrace |a|_pd, |b|_p\rbrace $. Now, note that for $p \notin S_{a,b}$, we have

\begin{align*}
\displaystyle  \int_{\overline{E}_p} \log^+ |t|_p + \log^+|at+b|_p d\mu_{\overline{E}_p}(t) = \sup_{|z|_p \leq 1} \log^+|z|_p + \sup_{|z|_p \leq 1} \log^+|az +b|_p = \log^+ \max \lbrace |a|_p, |b|_p \rbrace = 0.
\end{align*}

Therefore, there are finite many $v \in M_\mathbb{\mathbb{Q}}$ in (9). Furthermore, since 	$\mu_{\overline{E}_v} (\overline{E}_v) = \mu_{\overline{E}_v}(\mathbbm{1}_{\overline{E}_v})< \infty $ (where $\mathbbm{1}_{\overline{E}_v}$ is the continuous function $\mathbbm{1}_{\overline{E}_v}: \overline{E}_v \to \mathbb{R}$, given by $\mathbbm{1}_{\overline{E}_v}(t) = 1$ for each $t \in \overline{E}_v$) and $f(t) = \log^+|t|_v +\log^+|at+b| \leq \log^+r_i + \log^+ \max \lbrace |a|_pr_i, |b|_p \rbrace $,  for each  $ t \in \overline{E}_v$, we conclude that $M_{E_x} < \infty$.  Consequently, we have that

\begin{align*}
\mu^{ess}(h_{a,b}) &\leq \displaystyle \sum_{\nu \in M_\mathbb{\mathbb{Q}} } \int_{\overline{E}_\nu} \log^+ |t|_\nu + \log^+|at+b|_\nu d\nu_{\overline{E}_\nu}(t)\\
&= \displaystyle \sum_{i=1}^s  \sup_{|z|_{p_i} \leq r_i}\log^+ |z|_{p_i} +  \sup_{|z| \leq r_i} \log^+|az+b|_{p_i} + \frac{1}{2\pi} \int_0^{2\pi} \log^+ \left| \frac{e^{i\theta}}{r_1r_2... r_s} +x \right| + \log^+\left|\frac{ae^{i\theta}}{r_1r_2... r_s}+b+ax \right| d\theta.\\
&= \displaystyle \sum_{i=1}^s  \log^+ r_i +  \sup_{|z| \leq r_i} \log^+|az+b|_{p_i} + \frac{1}{2\pi} \int_0^{2\pi} \log^+ \left| \frac{e^{i\theta}}{r_1r_2... r_s} +x \right| + \log^+\left|\frac{ae^{i\theta}}{r_1r_2... r_s}+b+ax \right| d\theta\\
&=  \displaystyle \sum_{i=1}^s  \log^+ r_i +  \log^+ \max \lbrace |a|_{p_i}r_i,|b|_{p_i}\rbrace + \frac{1}{2\pi} \int_0^{2\pi} \log^+ \left| \frac{e^{i\theta}}{r_1r_2... r_s} +x \right| + \log^+\left|\frac{ae^{i\theta}}{r_1r_2... r_s}+b+ax \right| d\theta.\\
&= \Gamma_{a,b}(x,r_1,r_2,... , r_s).
\end{align*}

Moreover, suppose that $r_1,r_2,..., r_s \in \mathbb{Q}^+$ are fixed, and $j \in [\Gamma_{a,b}(x,r_1,r_2,..., r_s), \infty)$. Since,  $\Gamma_{a,b}$ is continuous and $\Gamma_{a,b}(x,r_1,r_2,..., r_s) \to \infty$, when $x \to \infty$, we can find $g \in \mathbb{R}$, such that, $j=\Gamma_{a,b}(g,r_1,r_2,...,r_s)= M_{E_g}$, and then, there is a sequence of algebraic numbers $(\alpha_{n})_{n \in \mathbb{N}}$, such that, $h_{a,b}(\alpha_{n}) \to j$. Hence, for each $x \in \mathbb{R}$, the image of $h_{a,b}$ is dense in the interval $[\Gamma_{a,b}(x,r_1,...,r_s), \infty)$. This concludes the proof of the theorem. $\hfill\square$\\

\hspace{1cm}

Experimental results show that the minimal value of $\Gamma$ is achieved when $r_1=r_2=...=r_s =1$, so, in general we will always take these values, furthermore, if we take $|a|=1$, then

\begin{center}
$\Gamma_{1,b} (x,1,1,...,1)= \Delta(1,b) + \displaystyle\frac{1}{2\pi} \int_0^{2\pi} \log^+ \left| e^{i\theta} +x \right| + \log^+\left|e^{i\theta}+b+x \right| d\theta = \Omega_{1,b}(x)$.
\end{center}

With this last observation and Theorem 4.2 we obtain Theorem C.  Using Theorem 4.2 together with Proposition 2.4 and Theorem 3.2, we can conclude Theorem A. Finally, using Proposition 2.8, Proposition 3.4 and Theorem 4.2, we obtain the following result \\ \\

\textit{{\bf{Theorem 4.3:}} Let $a,b \in \mathbb{Q}$, $|a| \geq 1$, $|b|-|a| > 1$, and $|b/a| \geq 4 $, then}

\begin{center}
$\displaystyle\log \left(  \frac{|b|-1}{|a|} \right) \leq \mu^{ess}(h_{a,b}) \leq \log \left( \frac{|b|}{|a|} \right) $.
\end{center}

\begin{center}
\bibliographystyle{plain}

\bibliography{BIBI}

\begin{thebibliography}{1}

\bibitem{MR2216774}
Enrico Bombieri and Walter Gubler.
\newblock {\em Heights in {D}iophantine geometry}, volume~4 of {\em New
  Mathematical Monographs}.
\newblock Cambridge University Press, Cambridge, 2006.

\bibitem{MR3802441}
Jos\'{e}~Ignacio Burgos~Gil, Ricardo Menares, and Juan Rivera-Letelier.
\newblock On the essential minimum of {F}altings' height.
\newblock {\em Math. Comp.}, 87(313):2425--2459, 2018.

\bibitem{MR3928039}
Jos\'{e}~Ignacio Burgos~Gil, Patrice Philippon, Juan Rivera-Letelier, and
  Mart\'{\i}n Sombra.
\newblock The distribution of {G}alois orbits of points of small height in
  toric varieties.
\newblock {\em Amer. J. Math.}, 141(2):309--381, 2019.

\bibitem{MR1838073}
Christophe Doche.
\newblock Zhang-{Z}agier heights of perturbed polynomials.
\newblock volume~13, pages 103--110. 2001.
\newblock 21st Journ\'{e}es Arithm\'{e}tiques (Rome, 2001).

\bibitem{riveraspectrum}
Shenxiong Li.
\newblock On the spectrum and essential minimum of heights in projective plane.
\newblock 2019.
\newblock
  http://www.shenxionghomepage.com/\textsc{R}esearch\%20\textsc{P}age/1.\%20\textsc{A}rithmetic\%20\textsc{D}ynamics/\textsc{O}n\%
  20the\%20\textsc{S}pectrum\%20and\%20\textsc{E}ssential\%20\textsc{M}inimum\%20of\%20\textsc{H}eights\%20in\%20\textsc{P}rojective\%
  20\textsc{P}lane.pdf.

\bibitem{MR1009368}
Robert~S. Rumely.
\newblock {\em Capacity theory on algebraic curves}, volume 1378 of {\em
  Lecture Notes in Mathematics}.
\newblock Springer-Verlag, Berlin, 1989.

\bibitem{MR1197513}
D.~Zagier.
\newblock Algebraic numbers close to both {$0$} and {$1$}.
\newblock {\em Math. Comp.}, 61(203):485--491, 1993.

\bibitem{MR1189866}
Shouwu Zhang.
\newblock Positive line bundles on arithmetic surfaces.
\newblock {\em Ann. of Math. (2)}, 136(3):569--587, 1992.

\end{thebibliography}
\end{center}

\end{document}